\documentclass[11pt,reqno]{amsart}

\usepackage{hyperref}
\usepackage{graphicx}
\usepackage{color}
\usepackage{amsfonts,amssymb}
\usepackage{bbm} 
\usepackage{psfrag}

\usepackage{mathrsfs}

\usepackage{a4wide}

\newtheorem{neu}{}[section]

\newtheorem*{Cor*}{Corollary}
\newtheorem{Thm}[neu]{Theorem}
\newtheorem*{Thm*}{Theorem}

\newtheorem*{Prop*}{Proposition}
\theoremstyle{definition}

\newtheorem*{Rmk*}{Remark}
\newtheorem{Rmk}[neu]{Remark}

\newtheorem*{Ex*}{Example}

\newtheorem*{Qu*}{Question}

\newcommand{\Z}{\mathbb{Z}}
\newcommand{\R}{\mathbb{R}}

\newcommand{\RP}{\R\mathrm{P}}

\newcommand{\pf}{\longrightarrow}

\newcommand{\CZ}{\mu_{\mathrm{CZ}}}

\newcommand{\om}{\omega}
\newcommand{\Om}{\Omega}

\newcommand{\A}{\mathcal{A}}

\renewcommand{\S}{\mathfrak{S}}

\newcommand{\D}{\mathbb{D}}

\renewcommand{\L}{\mathscr{L}}

\newcommand{\RFH}{\mathrm{RFH}}
\newcommand{\RFC}{\mathrm{RFC}}

\newcommand{\beq}{\begin{equation}}
\newcommand{\beqn}{\begin{equation}\nonumber}
\newcommand{\eeq}{\end{equation}}

\newcommand{\bea}{\begin{equation}\begin{aligned}}
\newcommand{\bean}{\begin{equation}\begin{aligned}\nonumber}
\newcommand{\eea}{\end{aligned}\end{equation}}

\numberwithin{equation}{section}

\definecolor{Urs}{rgb}{0,.7,0}
\definecolor{Peter}{rgb}{0,0,1}
\definecolor{red}{rgb}{1,0,0}

\newcommand{\p}{\partial}

\begin{document}
\title{On a Theorem by Ekeland-Hofer}
\author{Peter Albers}
\author{Urs Frauenfelder}
\address{
    Peter Albers\\
    Department of Mathematics\\
    Purdue University}
\email{palbers@math.purdue.edu}
\address{
    Urs Frauenfelder\\
    Department of Mathematics and Research Institute of Mathematics\\
    Seoul National University}
\email{frauenf@snu.ac.kr}
\begin{abstract}
In \cite[Theorem 1]{Ekeland_Hofer_Two_symplectic_fixed_point_theorems_with_applications_to_Hamiltonian_dynamics} Ekeland-Hofer prove that for a centrally symmetric, restricted contact type hypersurface in $\R^{2n}$ and for any global, centrally symmetric Hamiltonian perturbation there exists a leaf-wise intersection point. In this note we show that if we replace restricted contact type by star-shaped there exists infinitely many leaf-wise intersection points or a leaf-wise intersection point on a closed characteristic.
\end{abstract}
\maketitle

\section{Introduction}

Let $S\subset\R^{2n}$ be a hypersurface. Then $S$ carries a rank-1-foliation where the tangent space to a leaf $L_S(x)$ through $x\in S$ is given by  $\L_S(x):=\{v\in T_x S\mid \om(v,w)=0\;\;\forall w\in T_xS\}$. Here $\om$ is the standard symplectic form on $\R^{2n}$. A point $x\in S$ such that
\beq
\psi_1(x)\in L_S(x)
\eeq
is called a leaf-wise intersection point, see \cite{Moser_A_fixed_point_theorem_in_symplectic_geometry}. A hypersurface is of restricted contact type if there exists a 1-form $\lambda\in\Om^1(\R^{2n})$ with 
\beq
\begin{cases}
d\lambda=\om\\
\lambda_x(v)\neq0\quad\forall v\in\L_S(x)\;.
\end{cases}
\eeq
We call a hypersurface $\Z/2$-invariant or centrally symmetric if it is invariant under the symplectic involution $I:\R^{2n}\pf\R^{2n}$ given by $I(x)=-x$. We set $D_{\om,\Z/2}:=\{\phi\in\mathrm{Symp}(\R^{2n})\mid \phi\circ I=I\circ \phi\}$. In \cite{Ekeland_Hofer_Two_symplectic_fixed_point_theorems_with_applications_to_Hamiltonian_dynamics} Ekeland and Hofer prove the following theorem.

\begin{Thm}[\cite{Ekeland_Hofer_Two_symplectic_fixed_point_theorems_with_applications_to_Hamiltonian_dynamics}, Theorem 1]\label{thm:Ekeland_Hofer}
 Assume that $S\subset\R^{2n}$ is a connected, compact, $\Z/2$-invariant hypersurface of restricted contact type. Let $t\to\psi_t$ be an isotopy of the identity in $D_{\om,\Z/2}$. Then there exists a leafwise intersection point $x\in S$.
\end{Thm}

In this article we improve Theorem \ref{thm:Ekeland_Hofer} under the additional assumption that $S$ bounds a star-shaped (with respect to the origin) region in $\R^{2n}$, that is, it is of restricted contact type with respect to the standard primitive $\lambda_0=\tfrac12\sum x_idy_i-y_idx_i$.

\begin{Thm}\label{thm:new_Ekeland_Hofer}
 Assume that $S\subset\R^{2n}$ is a connected, compact, $\Z/2$-invariant, star-shaped hypersurface. Let $t\to\psi_t$ be an isotopy of the identity in $D_{\om,\Z/2}$. Then there exist infinitely many leaf-wise intersection points on $S$ or there exists a leaf-wise intersection point $y$ such that the leaf $L_S(y)$ is closed, that is, $L_S(y)$ is a closed characteristic.
\end{Thm}

\begin{Rmk}
If $n\geq2$ then for a generic isotopy of $\Z/2$-equivariant $\psi_t\in\D_{\om,\Z/2}$ there are no leaf-wise intersection points on closed characteristics. Hence, there exist infinitely many leaf-wise intersection points. This follows from a $\Z/2$-invariant version of \cite[Theorem 3.3]{Albers_Frauenfelder_Leafwise_Intersections_Are_Generically_Morse}. That Theorem 3.3 holds in the $\Z/2$-invariant case is due to the fact that for critical points $(v,\eta)$ of $\A$ (see below) the loop $v$ does not pass through the fix point $0$ of $I$ since for an invariant Hamiltonian function $H$ the Hamiltonian flow $\phi_H^t$ fixes $0$ for all times.
\end{Rmk}

\subsubsection*{Acknowledgments}
This article was written during visits of the authors at the Institute for Advanced Study, Princeton. The authors thank the Institute for Advanced Study for their stimulating working atmospheres. The authors are grateful to Helmut Hofer for many inspiring discussions. 

This material is based upon work supported by the National Science Foundation under agreement No.~DMS-0635607 and DMS-0903856. Any opinions, findings and conclusions or recommendations expressed in this material are those of the authors and do not necessarily reflect the views of the National Science Foundation.

\section{Equivariant Rabinowitz Floer homology}
We consider the standard symplectic space $(\R^{2n},\om=d\lambda_0)$ and the symplectic involution $I:\R^{2n}\pf\R^{2n}$ given by $I(x)=-x$. Let $F:\R^{2n}\pf\R$ be the $I$-invariant function $F(x):=\tfrac12\big(|x|^2-1\big)$. In particular, the Rabinowitz action functional
\bea
\A:C^{\infty}(S^1,\R^{2n})\times\R&\pf\R\\
(v,\eta)&\mapsto -\int v^*\lambda_0-\eta\int F(v)dt
\eea
is invariant under $I$. Moreover, $I$ acts freely on the critical points of $\A$ and on the space of gradient flow lines (in the sense of Floer) which asymptotically converge to critical points, where we use an $I$-invariant compatible almost complex structure $J$ to define the gradient of $\A$. Therefore, we can construct equivariant Rabinowitz Floer homology easily as follows:

\bea
\RFC^{\Z/2}_k(S^{2n-1},\R^{2n})&:=\RFC_k(S^{2n-1},\R^{2n})\Big/(\Z/2)\\
\p^{\Z/2}[x]&:=[\p x]
\eea
For details on the construction of Rabinowitz Floer homology and its relation with leaf-wise intersection points we refer to 
\cite{Cieliebak_Frauenfelder_Restrictions_to_displaceable_exact_contact_embeddings,Albers_Frauenfelder_Leafwise_intersections_and_RFH}
\begin{Thm}[\cite{Cieliebak_Frauenfelder_Restrictions_to_displaceable_exact_contact_embeddings}]\label{thm:RFH_0}
Since $S^{2n-1}$ is Hamiltonianly displaceable
\beq
\RFH_*(S^{2n-1},\R^{2n})\cong0
\eeq
\end{Thm}

\begin{Thm}
For all $k\in\Z$ we have
\beq
\RFH^{\Z/2}_k(S^{2n-1},\R^{2n})\cong\Z/2\;.
\eeq
\end{Thm}

\begin{proof}
For a critical point $(v,\eta)$ of $\A$ the $\eta$-periodic loop $v(t/\eta)$ is a Reeb orbit on $S^{2n-1}$ with respect to the standard contact form or in case $(v,0)$ the loop $v(t)$ is constant and represents a point on $S^{2n-1}$. Since the Reeb flow $\varphi^t$ on $S^{2n-1}$ is periodic the action functional $\A$ is Morse-Bott with critical manifolds
\beq
C_k\cong S^{2n-1}\quad k\in\Z
\eeq
where a point $x\in S^{2n-1}$ is identified with $(t\mapsto \varphi^{2\pi kt}(x),2\pi k)\in C_k$. The Conley-Zehnder index $\CZ$ equals $2nk$ on $C_k$. We fix on $S^{2n-1}$ the Morse function
\beq
f(x_1,\ldots,x_{2n}):=\sum_{i=1}^{2n}ix_i^2\;.
\eeq
$f$ descends to a $\Z/2$-perfect Morse function $\bar{f}$ on $\RP^{2n-1}=S^{2n-1}\big/(\Z/2)$. In particular, $\bar{f}$ has precisely one critical point in each degree $0,\ldots,2n-1$ and therefore, $f$ has critical points $y_l,z_l$ of degree $l=0,\ldots,2n-1$ satisfying $-y_l=z_l$. The Morse differential $\delta$ computes to
\beq
\delta y_l=y_{l-1}+z_{l-1}=\delta z_l,\quad\forall\;l=1,\ldots,2n-1
\eeq
and therefore in the quotient using the notation $\xi_l:=[y_l]=[z_l]$
\beq
\delta \xi_l=2\xi_{l-1}=0\quad\forall\;l=1,\ldots,2n-1\;.
\eeq
We define Morse functions
\beq
f^k:C_k\pf\R,\quad k\in\Z
\eeq
by $f^k:=f$ via the identification $C_k\cong S^{2n-1}$. We denote the critical points of $f^k$ by $y^k_l,z^k_l$, $l=0,\ldots,2n-1$.

The boundary operator $\p$ in Rabinowitz Floer homology is defined by counting gradient flows lines with cascades, see \cite{Frauenfelder_Arnold_Givental_Conjecture,Cieliebak_Frauenfelder_Restrictions_to_displaceable_exact_contact_embeddings}. Since the Conley-Zehdner index equals $2nk$ on the critical manifolds $C_k$ the complex $\RFC_*(S^{2n-1},\R^{2n})$ has exactly two generators in each degree. By index and energy reasons
\beq\label{eqn:floer_morse_diff_almost_always}
\p y_l^k=\delta y_l^k\text{ and }\p z_l^k=\delta z_l^k\quad\forall l=1,\ldots,2n-1,\;\forall k\in\Z\;.
\eeq
Again, by index reasons and by symmetry there exists $a^k,b^k\in\Z/2$ with
\beq
\p y_{0}^{k+1}=a^k y_{2n-1}^{k}+b^kz_{2n-1}^k=\p z_{0}^{k+1}\quad\forall k\in\Z\;.
\eeq
From $\p\circ\p=0$ and \eqref{eqn:floer_morse_diff_almost_always} we conclude $a^k=b^k$. According to Theorem \ref{thm:RFH_0} by Cieliebak-Frauenfelder we have $\RFH_*(S^{2n-1},\R^{2n})\cong0$. This implies that $a^k=b^k=1$ since otherwise $y_0^{k+1}$ is a cycle but not a boundary, compare figure \ref{fig:ladder}.

\begin{figure}[htb]
\psfrag{y1}{$y_l^k$}
\psfrag{y2}{$y_{l-1}^k$}
\psfrag{y3}{$y_{l-2}^k$}
\psfrag{z1}{$z_l^k$}
\psfrag{z2}{$z_{l-1}^k$}
\psfrag{z3}{$z_{l-2}^k$}
\psfrag{x1}{$\xi_l^k$}
\psfrag{x2}{$\xi_{l-1}^k$}
\psfrag{x3}{$\xi_{l-2}^k$}
\includegraphics[scale=1.5]{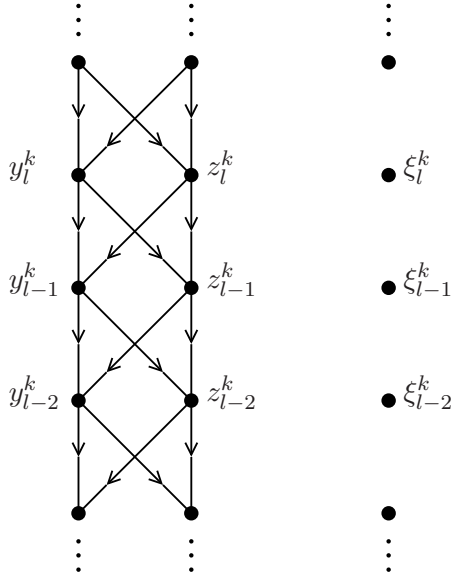}
\caption{The non-equivariant and the equivariant chain complexes.}\label{fig:ladder}
\end{figure}

With this we can compute the $\Z/2$-equivariant complex $\big(\RFC^{\Z/2}_k(S^{2n-1},\R^{2n}), \p^{\Z/2}\big)$ as follows. We have generators $\xi^k_l:=[y^k_l]=[z^k_l]$, $l=0,\dots,2n-1$, $k\in\Z$, of degree $\deg\xi_l^k=l+2nk$. In particular, there is exactly one critical point in each degree. We compute
\bea
\p^{\Z/2}\xi^k_l=[\p y^k_l]&=\begin{cases}
[y^k_{l-1}+z^k_{l-1}]&\text{for }l=1,\ldots,2n-1\\
[y_{2n-1}^{k-1}+z_{2n-1}^{k-1}]&\text{for }l=0
\end{cases}\\
&=\begin{cases}
2\xi^k_{l-1}&\text{for }l=1,\ldots,2n-1\\
2\xi^{k-1}_{2n-1}&\text{for }l=0
\end{cases}\\
&=0\;.
\eea
That is, the equivariant complex is acyclic. This proves the Theorem.
\end{proof}

In \cite{Albers_Frauenfelder_Spectral_invariants_in_RFH} we associated spectral values $\sigma(\xi)$ to homology classes $\xi$ in Rabinowitz Floer homology. We define
\beq
\S:=\{\sigma(\xi_l^k)\mid\xi_l^k\in\RFH_*^{\Z/2}(S^{2n-1},\R^{2n})\}\;.
\eeq
From the proof of Theorem \ref{thm:RFH_0} is follows immediately
\beq\label{eqn:spectrum_is_big}
\S=2\pi\Z
\eeq
since $\A(\xi_l^k)=-2\pi k$.

\section{Proof of Theorem \ref{thm:new_Ekeland_Hofer}}

We first assume that the isotopy $\psi_t$ is generated by a compactly supported Hamiltonian function $H:\R^{2n}\times [0,1]\pf\R$. Since $\psi_t\circ I=I\circ \psi_t$ we can assume
\beq
H(t,I(x))=H(t,x)\;.
\eeq
Moreover, since $S$ is star-shaped it is a graph over the standard sphere $S^{2n-1}$. Therefore, we can find a family of functions $F_r:\R^{2n}\pf\R$, $r\in[0,1]$, such that $F_1^{-1}(0)=S$, $F_0=F=\tfrac12\big(x^2-1\big)$, and all hypersurfaces $F_r^{-1}(0)$ are $I$-invariant and graphs over $S^{2n-1}$. Thus, all Rabinowitz action functionals
\bea
\A_r:C^{\infty}(S^1,\R^{2n})\times\R&\pf\R\\
(v,\eta)&\mapsto -\int v^*\lambda_0-\eta\int F_r(v)dt-\int rH(t,v)dt
\eea
are $I$-invariant. Moreover, $I$ acts freely on the critical points and gradient flow lines for each $r\in[0,1]$. Equality \ref{eqn:spectrum_is_big} implies that $\A_0$ has critical points of arbitrarily large critical value. \cite[Corollary 5.13]{Albers_Frauenfelder_Spectral_invariants_in_RFH} implies that
then also $\A_1$ has to have critical points with arbitrarily large critical value. In particular, $\A_1$ has infinitely many critical points. It follows from \cite[Proposition 2.4]{Albers_Frauenfelder_Leafwise_intersections_and_RFH} that critical points of $\A_1$ give rise to leaf-wise intersections. Moreover, the map from critical points to leaf-wise intersection points is injective unless there exists a leaf-wise intersection on a closed characteristic. This proves the theorem in case that $\psi_t$ is generated by a compactly supported Hamiltonian function.

A general isotopy $\psi_t\in\mathrm{Symp}(\R^{2n})$ is generated by a Hamiltonian function $\widetilde{H}:\R^{2n}\times[0,1]\pf\R$ which however is not necessarily compactly supported. The set
\beq
K:=\{\psi_t(x)\mid x\in S,\;t\in[0,1]\}\subset\R^{2n}
\eeq
is compact since $S$ is compact. In particular, all critical points $(v,\eta)$ of $\A$ satisfy
\beq
v(t)\in K\quad\forall t\in S^1\;.
\eeq
Thus, we can cut-off $\widetilde{H}$ outside $K$ to make it into a compactly supported Hamiltonian without changing the critical points of $\A$. Thus, by first part of the proof we are done.

%
%
\bibliographystyle{amsalpha}
\bibliography{../../../Bibtex/bibtex_paper_list}
\end{document}